\documentclass[a4paper,11pt]{article}

\usepackage{amscd}
\usepackage{amsfonts}
\usepackage{amssymb}
\usepackage{euscript}
\usepackage{amsmath}
\usepackage{amsthm}
\usepackage{dsfont}

\newcommand{\Oc}{\mathcal{O}}
\newcommand{\na}{\otimes^{L}}
\newcommand{\at}{\mathop{\rm at}}
\newcommand{\id}{{\rm id}}
\newcommand{\Ext}{\mathop{\underline{RHom}^\bullet}}
\newcommand{\D}{\mathrm D(X)}
\newcommand{\Tc}{\mathcal T}
\newcommand{\Uc}{\EuScript U}
\newcommand{\Fc}{\EuScript F}
\newcommand{\Mc}{\mathcal M}
\newcommand{\I}{{\bf I}}
\newcommand{\E}{{\bf E}}
\newcommand{\Lc}{{\bf L}}
\newcommand{\Dc}{\EuScript  D}
\newcommand{\Tr}{\mathop{\rm Tr}}
\newcommand{\Td}{\mathop{\rm Td}}
\newcommand{\Ch}{\mathop{\rm Ch}}
\newcommand{\td}{\mathop{\rm td}}
\newcommand{\ch}{\mathop{\rm ch}}
\newcommand{\jet}{\mathop{\rm J^1}}
\newcommand{\At}{\mathrm{At}}
\newcommand{\K}{{\mathbf K}}

\newtheorem{prop}{Proposition}
\newtheorem{theorem}{Theorem}
\newtheorem{lemma}{Lemma}
\newtheorem{definition}{Definition}

\title{The Atiyah class, Hochschild cohomology and the Riemann-Roch theorem}
\author{Nikita Markarian}

\begin{document}

\maketitle

\begin{abstract}
We develop a formalism involving Atiyah classes of sheaves on a smooth manifold, Hochschild chain and cochain complexes. As an application we prove a version of the Riemann-Roch theorem.
\end{abstract}

\section*{Introduction}
The present paper grew out  of the question posed to the author by
B.~Feigin: "Why does the Todd class look like the invariant volume
form on a Lie group?" Our answer is contained in the
proof of Proposition \ref{td}.

In Section \ref{alg diff calc} we develop a formalism describing Atiyah
classes, Hoch\-schild (co)homology and relationships between
them. We introduce a global analog of the Hochschild-Kostant-Rosenberg isomorphism (\cite{HKR}). This construction has appeared in the literature, see \cite{GS},
\cite{Sw} and \cite{Y}. Our definition of the global HKR isomorphism appeared in the preprint of the present paper and was used in \cite{Ca} and \cite{Ra}.

Let $\D$ denote the derived category of sheaves of $\Oc$-modules on a smooth variety $X$.
One can consider the Atiyah class as a morphism (in the category of endofunctors of $\D$) from the identity functor to the functor of tensoring by the cotangent bundle shifted by one:
\begin{equation}
\at\! \colon \, \id \to \,\cdot\,\otimes \Omega^1[1]
\label{0}
\end{equation}
(see \eqref{com}). One may think about \eqref{0} as an action of the object $T[-1]$ dual to $\Omega^1[1]$  on  the identity functor of $\D$.
Iterating this action one can make the tensor power of $T[-1]$ (in fact, the symmetric power) act on the identity functor of $\D$. This is an heuristic way to look at the map $\I$ in \eqref{I} which relates the Atiyah class with the Hochschild cochain complex.

Most of the constructions of Section \ref{alg diff calc} are applicable in a more
general situation. Let $D$ be a closed symmetric monoidal category with a compatible triangulation (see \cite{M}). That is, $D$ is a triangulated category with  a symmetric product given by a functor $\otimes\colon D\times D \to D$ which is exact in each variable (with some compatibility properties that we do not need).
In this setting we introduce the \emph{category of  K\"ahler differentials} of $D$ as follows. The objects of this category are pairs $(M,\alpha)$, where $M$ is an object of $D$ and $\alpha$ is a morphism from the identity functor to  the functor of tensoring by $M$:
\[
\alpha \colon\, \id \to \,\cdot\,\otimes M ,
\]
such that, for any $E$ and $F$ in $D$,
\[
\alpha(E\otimes F)=\alpha(E)\otimes \id + \id\otimes \alpha (F),
\]
(More generally, instead of the functor $\cdot\,\otimes M$, one can consider an endofunctor $\Mc$ of $D$ such that for any $E,F \in D$ $\Mc(E\otimes F)=E\otimes
\Mc(F) $.) Morphisms in the  category of K\"ahler differentials are morphisms between the $M$'s which respect the $\alpha$'s. Assume that the category of K\"ahler differentials has an initial object given by $(\Omega^1 [1], At)$ (for suitable $\Omega^1$ in $D$). In this case we will refer to the object $\Omega^1$ as the {\em cotangent complex} of $D$
and the morphism of functors $At$ as  the \emph{Atiyah class}.

Consider the morphism
\[
At\colon\, \Omega^1[1] \to \Omega^1\otimes \Omega^1[2] .
\]
One can show that it defines a structure of a Lie coalgebra (in $D$) on $\Omega[1]$.
Therefore, the dual object $T[-1]$ has a structure of a Lie algebra which acts on the identity functor of $D$.
By analogy with the usual Lie algebra theory one may define the enveloping algebra of $T[-1]$ (an associative algebra in $D$) and its action extending the action of the $T[-1]$.

The basic example of such a situation is the derived category  of representations
of a Lie algebra $\mathfrak g$. In this case, we have $\Omega^1[1] = \mathfrak{g}^\vee$, i.~e. the (shifted) cotangent complex is the coadjoint representation $\mathfrak{g}^\vee$. For a representation $V$ the
Atiyah class is given by the map
\[
V\to V\otimes {\mathfrak g}^\vee
\]
induced by the action of $\mathfrak g$ on $V$. 
The dual object is $\mathfrak{g}$ itself
and the Lie algebra structure on it is the usual one. Its
enveloping algebra is the universal enveloping algebra of $\mathfrak g$
in the usual sense equipped with the adjoint action.

Another example is the subject of the first part of the paper.
The category $D$ is the derived category of sheaves of $\Oc$-modules on a smooth variety
with the usual tensor product. Comparing with the basic example one sees that the Hochschild cochain complex corresponds to the universal enveloping algebra,
the Hochschild chain complex corresponds to formal functions
on the group in the neighborhood of the unit and so on. The reader may find this analogy helpful.

In the second part of the paper we prove the Riemann-Roch theorem as an application the techniques developed in the first part. 

Essentially, we follow  \cite{TT}. However, instead of explicit
calculations with the \v{C}ech cocycles we work in the derived
category and use our algebraic-differential calculus.
The proof consists of two parts. In first part we reduce the theorem to a calculation of the dual class of the diagonal in the terminology of \cite{TT}. In the second part we perform the calculation.

\smallskip
{\bf Acknowledgments.}
It is а pleasure to express my gratitude to many people. It was B.~Feigin who posed the question and discussed it a lot. I thank M.~Kontsevich for many conversation about Hochschild cohomology. Constant interest  of A.~Tyurin at the early stages of this work was very encouraging. I am grateful to the participants of the seminar at the University of Chicago for their interest and especially to V.~Drinfeld who made a crucial remark. I am grateful to P.~Bressler and T.~Bridgeland for their inspiring interest and friendly support and to S.~Willerton  for valuable remarks. I am grateful to the referee for significant improvement of the exposition and pointing out references.

\subsection*{Notations}
\noindent
 $X$ is a  smooth algebraic variety over a field $k$ of
characteristic $0$ or bigger than $\dim X$. Everything works  for the complex-analytic case as well.

\smallskip \noindent
By $\D$ we denote the derived category of perfect complexes of
sheaves of $\Oc$-modules on $X$.

\smallskip \noindent
$\Delta$ always means the diagonal and the diagonal embedding.

\smallskip \noindent $p_i $ means projection onto the $i$-th factor of
$X^n$.

\smallskip \noindent
For $E\in\D $ by $E_\Delta$ we denote $\Delta_* E=p^*_1E\otimes
\Oc_\Delta=p^*_2E\otimes \Oc_\Delta$, where $\Delta$ is the diagonal in $X\times X$.
% To emphasize that we consider the object as an object on
% $X$ not $X\times X$ we shall write $Rp_{1*} E_\Delta$.

\smallskip \noindent
For $E\in \D$ by  $E^\vee$ we denote the dual object $\Ext (E,\Oc)$.

\smallskip \noindent
By $\Tr$ we denote the canonical morphism $ E\na E^\vee\to \Oc$.

\smallskip \noindent
By $\Omega^i$  we denote the  bundle of exterior forms. $\Lambda^i T$
is dual to $\Omega^i$.

\smallskip \noindent
$\omega $ denotes the bundle of exterior forms of top degree.

\section{Algebraic-differential calculus}\label{alg diff calc}

\subsection{The Atiyah class}

Let  $\Delta\subset X\times X$ be the diagonal and let  $I$ denote the ideal
sheaf of $\Delta$. Then, by definition, $\Oc_\Delta=\Oc_{X\times X}/I$ and
$\Omega_\Delta=I/I^2$. The two-step filtration on $\Oc_{X\times X}/I^2$
by powers of $I$ gives rise to the exact sequence 
\begin{equation}
\begin{CD}
0 @>>> \Omega^1_\Delta @>>> \Oc_{X\times X}/I^2 @>>> \Oc_\Delta @>>> 0
\end{CD}
\label{d}
\end{equation}
Since the terms of the sequence \eqref{d} are supported on the diagonal one may consider \eqref{d} as a sequence of sheaves of $\Oc_X$-$\Oc_X$-bimodules on $X$.
The two $\Oc_X$-module structures coincide on $\Omega^1_\Delta$ and $\Oc_X$ but are different on the middle term.

Let $E$ be a sheaf of $\Oc$-modules or a complex of such sheaves on $X$.
Take its tensor product with \eqref{d} with respect to the left
$\Oc$ module structure and consider it as a right $\Oc$ module. In
other words, tensor \eqref{d} by $p_1^*E$ and take the direct image
$p_{2*}$. Because all terms in \eqref{d} are locally free left
$\Oc$-modules this operation is exact and one gets an exact sequence
\begin{equation}
\begin{CD}
0@>>> E\otimes\Omega^1@>>>\jet(E)@>>> E@>>>0.
\end{CD}
\label{d2}
\end{equation}
Here $\jet(E)$ denotes $E\otimes\Oc/I^2$ with the right $\Oc$-module
structure and is called the \emph{sheaf of the first jets}.

\begin{definition}[\cite{A,I}]
For a sheaf of $\Oc$-modules or a complex of such sheaves$E$ on $X$
the class of extensions represented by \eqref{d2} is called the \emph{Atiyah class } $\at(E)\in Ext^1(E, E\otimes\Omega^1)$ of $E$.
\end{definition}
The Atiyah class is the only obstruction to the existence of a connection on a sheaf (\cite{A}).

\begin{definition}
A \emph{connection} on a  sheaf of $\Oc$-modules $E$ is a splitting
of the exact sequence \eqref{d2}, i.e. a map
$\nabla\colon E\to \jet(E)$ whose composition with the projection
$\jet(E)\to E$  equals to the identity map.
\end{definition}
% If  $E$ is a vector bundle
% then the central term in (\ref{a}) is just bundle of first jets of $E$ and a splitting of (\ref{a})
% is by very definition a connection on $E$.

The Atiyah class of an object in the derived category is defined in a way compatible with the definition above. Let $\At$ be the morphism in the derived category represented by
extension \eqref{d}:
\begin{equation}
\At \colon \Oc_\Delta \longrightarrow \Omega^1_\Delta[1].
\end{equation}

\begin{definition}
For $E\in\D$ the morphism in $\D$
\[
\at (E)\colon E \to E\na\Omega^1[1]
\]
given by $Rp_{2*}(\At\na p_1^* E)$ is the {\em Atiyah class} of $E$.
\end{definition}

A trivial but important observation is that the Atiyah class is natural. That is,
for any morphism $E\stackrel{m}{\to} F$, the diagram
\begin{equation}
\begin{CD}
E            @>m>>          F \\
@V\at(E)VV               @V\at(F)VV \\
E\na\Omega^1[1] @>m\otimes\id>> F\na\Omega^1[1]
\end{CD}
\label{com}
\end{equation}
is commutative.

Here are two useful lemmas.

\begin{lemma}
Let $E=(E^i, d^i\colon E^i\to E^{i+1})$ be a complex of sheaves of $\Oc$-modules.
Asume given a connection $\nabla^i$ on $E^i$.
Then, $\at(E)$ is represented by
\[(\nabla d)^i \stackrel{def}{=} (d^i\circ\nabla^i-\nabla^{i+1}\circ d^i)\colon\, E^i\to E^{i+1}\otimes\Omega^1.\]
\label{con}
\end{lemma}

\begin{proof} It suffices to show that the Atiyah class of $E$ is given by the extension
\begin{equation}
\begin{CD}
0@>>> E\otimes\Omega^1 @>>> cone \,(\nabla d) @>>> E@>>> 0 \ .
\end{CD}
\label{cone}
\end{equation}
Consider the short exact sequence of complexes
\begin{equation}
\begin{CD}
   @. \dots   @. \dots @. \dots @. \\
@.     @AAA      @AAA     @AAA    @.\\
0 @>>>E^{i+1}\otimes\Omega^1  @>>> \jet(E^{i+1}) @>>> E^{i+1} @>>>0\\
@. @Ad^i\otimes\id AA       @A\jet(d^i)AA      @A d^iAA  @. \\
0 @>>> E^{i}\otimes\Omega^1@>>>\jet(E^{i}) @>>> E^{i}@>>>0\\
@.     @AAA      @AAA     @AAA   @.\\
 @. \dots   @. \dots @. \dots @.
\end{CD}
\label{big}
\end{equation}
which by definition represents the Atiyah class.
% class
% The central term in the exact sequence (\ref{a}) is isomorphic to
% $(j^1E^i, j^1d^i)$, where $j^1$ means first jets.
The connections on the terms of the complex $E$ give rise to the isomorphisms (splittings) $\jet(E^i)=E^i\oplus (E^i\otimes\Omega^1)$.
In terms of these splittings we have:
\[\jet(d^i)=
\begin{pmatrix} d^i& (\nabla d)^i
 \\ 0& d^i\otimes \id\end{pmatrix}\]
which means that \eqref{big}
is equal to \eqref{cone}.
This proves the lemma.
\end{proof}

\begin{lemma}  Let $E$ and $F$ be objects of
$\D$. Then $\at(E\na F)=\at(E)\otimes \id+\id\otimes\at(F)$.
\label{na}
\end{lemma}

\begin{proof} Let  $\Delta$ and $\Delta_i$ ($i\in \{1,2,3\}$)
be the submanifolds of $X^3$ defined by the equations
$x_1=x_2=x_3$ and $x_j=x_k$ for $j, k \ne i$ respectively, where $(x_1, x_2, x_3)$ is a point of $X^3$.
Let $I$ and $I_i$ denote the ideal sheaves  of  $\Delta$ and $\Delta_i$ in $X^3$.

The sheaf $I/I^2$ is by definition the conormal bundle to $\Delta$ in $X^3$. Therefore, it is isomorphic to the direct sum of two copies of $\Omega_\Delta^1$.
The sheaf $I/(I^2+I_i)$ is the restriction of the conormal bundle to $\Delta_i$.
Hence,  it is isomorphic to $\Omega^1_\Delta$.
There are three projections:
$\pi_i \colon I/I^2 \to I/(I^2+ I_i)=\Omega_\Delta^1$.
Identify $I/I^2$ with the direct sum of two copies of $\Omega_\Delta^1$
so that projections on  the first and the second summands are equal to
$\pi_1$ and $\pi_2$. Then, $\pi_3$ is equal to $-\pi_1-\pi_2$.

Consider the exact sequences
\begin{equation}
\begin{CD}
0 @>>> I/(I^2+I_i)@>>>\Oc/(I^2+I_i) @>>> \Oc_\Delta @>>> 0
\end{CD}
\end{equation}
and the corresponding morphisms in the derived category
\begin{equation}
\begin{CD}
\alpha_i\colon \Oc_\Delta @>>> I/(I^2+ I_i) [1]=\Omega^1_\Delta[1].
\end{CD}
\end{equation}
By  definition $Rp_{3*}(\alpha_i\na( p_1^* E\na p_2^*F))$
is equal to $\at (E)\na \id$, $\id \na \at (F)$ and $\at (E\na F)$
for $i=1,2,3$ respectively.

Consider the commutative diagram
\begin{equation}
\begin{CD}
0 @>>> I/I^2 @>>>     \Oc/I^2 @>>>    \Oc_\Delta @>>> 0\\
@.   @V\pi_iVV          @VVV          @V|VV          @.\\
0 @>>> I/(I^2+I_i)@>>>\Oc/(I^2+I_i) @>>> \Oc_\Delta @>>> 0
\end{CD}
\label{pr0}
\end{equation}
where the vertical arrows are the natural projections.
It gives a commutative diagram in the derived category
\begin{equation}
\begin{CD}
\Oc_\Delta @>>> I/I^2 [1] \\
@|            @V\pi_iVV \\
\Oc_\Delta @>\alpha_i>> I/(I^2+ I_i) [1]
\end{CD}
\label{pr}
\end{equation}
Applying $Rp_{3*}(-\na( p_1^* E\na p_2^*F))$
 to  the top row of \eqref{pr} we get the morphism
\begin{equation}
E\na F \to E\na F \na (\Omega^1\oplus\Omega^1)[1],
\label{ext}
\end{equation}
where $I/I^2$ is identified with $\Omega^1_\Delta\oplus\Omega^1_\Delta$ as above.
We denote the projections of  $\Omega^1\oplus\Omega^1$  onto the summands by $\pi_1$ and $\pi_2$ as well.
The compositions of  the extension \eqref{ext} with the projections $\pi_1$, $\pi_2$
and $\pi_3=-\pi_1-\pi_2$ are equal to $Rp_{3*}(\alpha_i\na (p_1^* E\na p_2^*F))$.
The latter are equal to $\at (E)\na \id$, $\id \na \at (F)$  and $\at (E\na F)$
for $i=1,2,3$ respectively.
This proves the lemma.

\end{proof}

\subsection{The Lie algebra $\Tc$}

Consider the Atiyah class of the cotangent bundle
\begin{equation}
\at(\Omega^1) \colon\Omega^1 \to \Omega^1\otimes \Omega^1[1]
\label{lie}
\end{equation}

\begin{prop}
\begin{enumerate}
\item $\at(\Omega^1)$ is symmetric, i.~e. invariant under the permutation of factors in
$\Omega^1\otimes\Omega^1$.
\item $\at(\Omega^1)$ obeys the \emph{Jacobi identity}, i.~e. the projection of
$\at(\Omega^1)\otimes\id \circ \at(\Omega^1)$
onto the part of $\Omega^1\otimes\Omega^1\otimes\Omega^1$ invariant under permutations is equal to zero.
\end{enumerate}
\end{prop}

\begin{proof} Let $\Delta$ be the diagonal in $X\times X$ and  $I$ be its ideal sheaf.

{\bf 1.} By definition,  the Atiyah class of $\Omega^1$  is given by
the image under $Rp_{2*}$ of the exact sequence
\begin{equation}
\begin{CD}
0 @>>> I/I^2\otimes p_1^*\Omega^1 @>>>  \Oc/I^2\otimes p_1^*\Omega^1 @>>>  \Oc /I\otimes p_1^*\Omega^1 @>>> 0
\end{CD}
\label{1}
\end{equation}
It is enough to show that this extension may be represented as the composition of some extension and the embedding of
 $(S^2\Omega^1)_\Delta$ into $(\Omega^1\otimes\Omega^1)_\Delta=I/I^2\otimes p_1^*\Omega^1$
as the symmetric part.

The sequence \eqref{1} is included in the commutative diagram
\begin{equation}
\begin{CD}
0 @>>> I^2/I^3 @>>> I/I^3 @>>> I/I^2 @>>> 0\\
@.      @VdVV        @VdVV     @VdV|V  @. \\ 
0 @>>> I/I^2\otimes p_1^*\Omega^1 @>>>  \Oc/I^2\otimes p_1^*\Omega^1 @>>>  \Oc/I\otimes p_1^*\Omega^1 @>>> 0
\end{CD}
\label{2}
\end{equation}
where vertical arrows are given by the exterior differential (\cite{EGA}).
It follows that the extension given by the bottom row
is composition of the one given by the top row and
the exterior differential $I^2/I^3\to I/I^2\otimes p_1^*\Omega^1$.
The latter composition is equal to the embedding of $(S^2\Omega^1)_\Delta=I^2/I^3$
in $(\Omega^1\otimes\Omega^1)_\Delta$ as the symmetric part.
Thus, the Atiyah class is symmetric.

{\bf 2.} In addition to \eqref{2} consider the commutative diagram
\begin{equation}
\begin{CD}
0 @>>> I^3/I^4 @>>> I^2/I^3 @>>> I^2/I^3 @>>> 0\\
@.      @VdVV           @VdVV            @VdV|V      @.            \\
0 @>>> S^2_\Delta\Omega^1\otimes p_1^*\Omega^1
@>>>  I/I^3\otimes p_1^*\Omega^1 @>>>  \Omega^1_\Delta\otimes p_1^*\Omega^1 @>>> 0
\end{CD}
\label{3}
\end{equation}
where vertical arrows are again given by the exterior differential
(\cite{EGA}), the left vertical map being the embedding of the symmetric tensor
product into the partially symmetrized tensor product. By the previous paragraph, the image under $Rp_{2*}$ of the bottom row gives the extension $Sym(\at
(\Omega^1)\otimes\id)$, where $Sym$ denotes the projection onto the
symmetric part of $\Omega^1\otimes\Omega^1$.

Therefore, to prove the lemma it is enough to show that push-out of  the composition
of extensions given by the bottom rows of \eqref{2} and \eqref{3} 
by the projection onto the symmetric part of
$S^2\Omega^1_\Delta\otimes p_1^*\Omega^1 $ is trivial.

Combining \eqref{2} and \eqref{3} one obtains the following commutative diagram in the derived category of $X^2$: 
 \begin{equation}
\begin{CD}
I/I^2 @>>>   I^2/I^3[1] @>>>  I^3/I^4[2] @=S^3\Omega^1_\Delta[2] \\
  @VdV|V           @VdVV            @VdVV  @| \\
\Oc_\Delta\otimes p_1^* \Omega^1 @>>>  \Omega^1\otimes p_1^*\Omega^1[1] @>>>  S^2\Omega^1_\Delta\otimes p_1^*\Omega^1 [2] @>>> S^3\Omega^1_\Delta[2]
\end{CD}
\label{4}
\end{equation}
Here, the last arrow in the bottom row is the projection onto the
symmetric part. Therefore, the composition of morphisms in the top row is
equal to the composition of morphisms in the bottom row. But
composition of first two arrows in the top row is zero because they
represent the extension classes between successive quotients
of  the two-step filtration on $I/I^4 $ given by
powers of $I$. This proves the lemma.
\end{proof}

It follows from the above proposition that there is a structure of a Lie (super)algebra in the derived category $\D$  
on the shifted tangent bundle $T[-1]$. The map dual
to \eqref{lie} is the bracket. Denote this Lie algebra by $\Tc$:
 \[
 [\ ,\ ] \colon \Tc\na\Tc \to \Tc.
 \]

 The following definition  could be given in more general context,
 but we need it only for a special case.

\begin{definition}
We say that an algebra $U$ in
$\D $
with unit $e\in Hom(\Oc, U) $ and multiplication
$m\colon U\na U \to U $ is the \emph{enveloping algebra} of $\Tc$ if the following holds
\begin{enumerate}
\item  There is a map $\imath \colon \Tc \to U$ such that
\[
\imath \circ [\ ,\ ]= (m - m^\sigma)\circ \imath\otimes \imath,
\]
where $m^\sigma$ is multiplication in the reverse order with an appropriate sign.
\item The map  defined as the composition of embedding
\[S^*\Tc\stackrel{def}{=} \bigoplus_i\Lambda^iT[-i]\hookrightarrow \bigoplus_i\Tc^{\otimes i}
\]
and multiplication
\[
\Tc\otimes\cdots\otimes\Tc\stackrel{\imath\otimes\cdots\otimes\imath}{\longrightarrow}
U\otimes\cdots\otimes U \rightarrow U
\]
is an isomorphism.
\end{enumerate}
\label{u}
\end{definition}

Note that the  enveloping algebra exists and is unique. The second
condition of Definition \ref{u} gives the underlying object of $\D$
and the first condition defines the multiplication on it. In fact, the
multiplication is given by the Campbell-Hausdorff  formula
\cite{BL}, which describes the co-product on the formal function ring in
exponential coordinates which is dual to the product on the universal
enveloping algebra in terms of the Poincar\'e-Birkhoff-Witt
isomorphism.

\subsection{Hochschild cohomology}
We define the object $\Uc$ in $\D$ (\emph{Hochschild cochain complex}) by
\begin{equation*}
\Uc \stackrel{def}{=} Rp_{1*} \Ext\nolimits_{X\times X}(\Oc_\Delta, \Oc_\Delta),
\end{equation*}
where $\Delta $ is the diagonal in $X\times X$.
$\Uc$ is endowed with the canonical structure of an algebra in $\D$.

Let $\pi\colon \Omega^1_{X\times X} \to p_1^* \Omega^1_X $ be the
natural projection. Let
\begin{equation}
\imath\colon \Tc \to \Uc
\label{imath}
\end{equation}
denote the map defined as the contraction with
\begin{equation}
(\id\otimes\pi)\circ \at(\Oc_\Delta)\in Ext^1(\Oc_\Delta, \Oc_\Delta\otimes p_1^* \Omega^1).
\label{ac}
\end{equation}

\begin{theorem}
The algebra  $\Uc$ together with the map $\imath $ is the  enveloping
algebra of $ \Tc$ in the sense of Definition \ref{u}.
% That is
% \begin{enumerate}
% \item
% \item
% \end{enumerate}
\label{un}
\end{theorem}

\begin{proof}

We prove that the first condition of Definition \ref{u} is satisfied.
To simplify notations we omit the projection $\pi$.

One needs to show that the symmetric part of the composition
\[
\Oc_\Delta \stackrel{\at(\Oc_\Delta)}{\longrightarrow} \Oc_\Delta\otimes p_1^* \Omega^1  [1]
\stackrel{\at(\Oc_\Delta)\otimes \id}{\longrightarrow}
  \Oc_\Delta\otimes p_1^* \Omega^1 \otimes p_1^*\Omega^1  [2]
\]
is equal to the composition
\[
\Oc_\Delta \stackrel{\at(\Oc_\Delta)}{\longrightarrow}
\Oc_\Delta\otimes p_1^* \Omega^1 [1]
\stackrel{\id\otimes \at(\Omega^1)}{\longrightarrow}
\Oc_\Delta\otimes p_1^* \Omega^1 \otimes p_1^*\Omega^1  [2]
\]
Consider the diagram
\[
\begin{CD}
\Oc_\Delta @>\at(\Oc_\Delta)>> \Oc_\Delta\otimes p_1^*\Omega^1 [1]\\
@V\at(\Oc_\Delta)VV                    @V\at(\Oc_\Delta\otimes p_1^*\Omega^1)VV \\
\Oc_\Delta\otimes p_1^*\Omega^1[1] @>\at(\Oc_\Delta)\otimes\id>> (\Oc_\Delta\otimes p_1^*\Omega^1) \otimes p_1^*\Omega^1  [2]
\end{CD}
\]
It is (super)commutative by \eqref{com}. By Lemma \ref{na} we have the equality
$\at(\Oc_\Delta\otimes p_1^*\Omega^1)=\at(\Oc_\Delta)\otimes\id+\id\otimes\at(p_1^*\Omega^1)$.
This proves the first condition.

To prove the second condition it suffices to show that the morphism from the second part of
Definition \ref{u} is an isomorphism at any geometric point of  $X$.
Because the construction of $\Uc$ is natural with respect to open embeddings it is sufficient to prove it for the spectrum of a local ring. Moreover, one can
replace the local ring with its completion, because completion is
flat. Thus, one may assume that $X=\mathop{\rm Spec} k[[x_1, \dots, x_n]]$. But in this case it
could be solved by direct calculation making use of  the Koszul
resolution. Also, Proposition \ref{hkr} below could be applied.
\end{proof}

We denote the isomorphism from the second part of Definition \ref{u} by
\begin{equation}
\I\colon  \bigoplus_i\Lambda^iT[-i] \stackrel{\sim}{\rightarrow} \Uc
\label{I}
\end{equation}

In \cite{HKR} the case of a smooth affine manifold $X=\mathop{\rm
Spec} (A)$ was considered. There  the cohomology of
 $\Uc$, which is to say, $\mathop{Ext}^i_{A\otimes A}(A, A) $ by means of the \emph{standard
 resolution} \cite{BH}.

The standard resolution $B=(B_i, d_i)$, $i\geqslant0$ of $A$ as a $A\otimes A$-module is as follows:
$B_n$ is a free $A\otimes A$-module generated by tensor power $A^{\otimes n}$ over the base field
and the  differential is given by
\begin{multline}
d(a_1\otimes a_2\otimes\cdots\otimes a_{n-1}\otimes a_n)= \\a_1(a_2\otimes\cdots\otimes a_n)-(a_1 a_2\otimes\cdots\otimes a_n)+
\dots+(-1)^n(a_1\otimes a_2\otimes\cdots\otimes a_{n-1}) a_n .
\end{multline}

The Hochschild-Kostant-Rosenberg isomorphism (\cite{HKR})
is given by
\begin{equation}
\begin{split}
\Lambda^i T & \to Hom_{A\otimes A}(S_i, A)=Ext^i_{A\otimes A}(A, A) \\
\partial_1 \wedge \dots\wedge \partial_i\! & : \quad
a_1\otimes\cdots\otimes a_i
\mapsto \sum_{\sigma\in \Sigma_i} (-1)^{sign\, \sigma}
\partial_{\sigma(1)} a_1\cdots \partial_{\sigma(i)} a_i
\end{split}
\label{HKR}
\end{equation}

\begin{prop}
For a smooth affine manifold $X=\mathop{\rm Spec} (A)$  the isomorphism
$\I$ coincides with \eqref{HKR}. \label{hkr}
\end{prop}

\begin{proof} The standard resolution is a resolution of the structure
sheaf of the diagonal $\Oc_\Delta$. The terms of the standard resolution
are free modules, hence they are equipped with canonical connections.
Applying Lemma \ref{con}  to the standard resolution one obtains the following expression for the Atiyah class of $\Oc_\Delta$:
\begin{equation*}
\begin{split}
 \at(\Oc_\Delta) \colon &(a_1\otimes a_2\otimes\cdots\otimes a_{n-1}\otimes a_n)
 \mapsto \\
 &da_1 (a_2\otimes\cdots\otimes a_{n-1}\otimes a_n)
 +(-1)^n(a_1\otimes a_2\otimes\cdots\otimes a_{n-1})da_n ,
\end{split}
\end{equation*}
where $d$ is the exterior differential. To finish the proof substitute the above formula into \eqref{ac}.
\end{proof}

\subsection{Hochschild homology}
Let $\Fc$ denote the object in $\D$  ({\em Hochschild chain complex}) defined by
\begin{equation}
\Fc \stackrel{def}{=}R p_{1*} (\Oc_\Delta \na \Oc_\Delta),
\label{tc}
\end{equation}
where $\Delta $ is the diagonal in $X\times X$

The canonical action of $\Uc$ on $\Fc$
\begin{equation}
\Dc\colon  \Uc\na \Fc\to \Fc
\label{action}
\end{equation}
 is given by
\begin{equation*}
\Uc=\Ext(\Oc_\Delta, \Oc_\Delta)\stackrel{-\otimes \id}{\longrightarrow}
\Ext(\Oc_\Delta\na \Oc_\Delta, \Oc_\Delta\na\Oc_\Delta)=\Ext(\Fc, \Fc).
\end{equation*}
Note that the composition of this morphism with $\imath$
\eqref{imath} gives an action of $\Tc$ on $\Fc$
\begin{equation}
\Tc\na Rp_{1*}(\Oc_\Delta\na\Oc_\Delta) \stackrel{\Dc}{\rightarrow} Rp_{1*}(\Oc_\Delta\na\Oc_\Delta)
 \label{actt}
\end{equation}
which is equal to \eqref{ac} tensored by $\Oc_\Delta$.

The canonical morphism $\Oc_\Delta \na \Oc_\Delta\to \Oc_\Delta \otimes \Oc_\Delta = \Oc_\Delta$
gives rise to the morphism
\begin{equation}
\varepsilon \colon \Fc \to \Oc_\Delta
\label{eps}
\end{equation}

\begin{prop}
The composition of the maps $\Dc$  and $\varepsilon$ defines a perfect pairing
\begin{equation}
\Uc\na \Fc\to \Oc.
\label{par}
\end{equation}
\end{prop}
\begin{proof}  The statement is local and may be proved by the same
considerations as the second part of the proof of Theorem \ref{un}.  \end{proof}

It follows from the proposition that $\Fc$ is dual to $\Uc$.
Denote the isomorphism dual to $\I$ by
\[\E\colon(\bigoplus_i\Lambda^iT[-i])^\vee =\bigoplus_i\Omega^i[i] \xleftarrow{\sim} \Fc \].

Let $L^n\in Hom(\bigoplus_i\Omega^i[i] , (\bigoplus_i\Omega^i[i])\otimes\Omega^1 [1])$
 denote the morphism defined as follows: its restriction to $\Oc$ is zero and its restriction to 
$\Omega^1$ is given by the composition
\[
L^n\colon\Omega^1\stackrel{\at\nolimits^n(\Omega^1)}{\longrightarrow}
\Omega^1\otimes {(\Omega^1)}^{\otimes i}[i]
\stackrel{\id\otimes\wedge}{\longrightarrow} \Omega^1\otimes\Omega^n [n]=\Omega^n\otimes\Omega^1[n]
\]
(compare \eqref{ae}); $L^n $ is extended to all of
$\bigoplus_i\Omega^i[i]$ by the Leibniz rule with respect to the wedge product.

Let
\[
\Lc \colon \quad    \bigoplus_i\Omega^i[i] \to (\bigoplus_i\Omega^i[i])\otimes\Omega^1 [1]
\]
denote the morphism defined by the formula
\begin{equation}
\Lc=\sum l_n  L^n
\label{formula}
\end{equation}
where $\sum l_n z^n=z/(e^z-1) $.

The following theorem provides a description of the action of $\Tc\subset\Uc$ on $\Fc$
which allows us to obtain the action of all of $\Uc$.

\begin{theorem}
The diagram
\[
\begin{CD}
\Tc\otimes\Fc @>\Dc>> \Fc \\
@V\id\otimes\E VV @V\E VV \\
T\otimes\bigoplus_i\Omega^i[i] [-1]@>\Lc>> \bigoplus_i\Omega^i[i]
\end{CD}
\]
is commutative. \label{L}
\end{theorem}

\begin{proof} The map adjoint to the action
of $\Tc$ on $\Fc$ with respect to \eqref{par} is simply the (right)
multiplication on $\Uc$. Thus, the question is reduced to the problem of describing the multiplication in terms of the isomorphism $\I$. This is a purely combinatorial question and the answer may be found in any handbook on Lie algebras (e.~g. \cite {BL}).
% [ Ch2 \P6 ex. 3 ]
In dual terms, the problem is to write down
the left invariant fields in exponential coordinates. \end{proof}

In \cite{HKR}, in the case of a smooth affine variety $X=\mathop{\rm Spec} (A)$, the isomorphism dual (with respect to \eqref{par}) to \eqref{HKR} was computed in terms of the standard resolution $S$:
\begin{equation}
\begin{split}
 S_i\otimes_{A\otimes A} A=Tor^i_{A\otimes A}(A, A)& \to\Omega^i  \\
a_0 (a_1\otimes\cdots\otimes a_i) &\mapsto a_0 \,
da_1\wedge\dots\wedge da_i
\end{split}
\label{mu}
\end{equation}

\begin{prop}
For a smooth affine manifold $X=\mathop{\rm Spec}
 (A)$ the isomorphism $\E$ coincides with the isomorphism
\eqref{mu}.
\end{prop}

\begin{proof} Follows from the proof of Proposition \ref{hkr} and
the definition of $\E$.
\end{proof}

There is a multiplication on $\Fc$ defined by the composition
\[
(\Oc_\Delta\na\Oc_\Delta)\na(\Oc_\Delta\na\Oc_\Delta) \stackrel{\id\otimes\sigma\otimes\id}{\longrightarrow} \Oc_\Delta\na\Oc_\Delta\na\Oc_\Delta\na\Oc_\Delta \stackrel{\varepsilon\otimes\varepsilon}{\longrightarrow} \Oc_\Delta\na\Oc_\Delta,
\]
where $\sigma$ is the permutation of factors. One can show that it corresponds to the usual multiplication on differential forms under isomorphism $E$.
We do not need this fact.

\section{The Riemann-Roch theorem}\label{RR thm}

\subsection{Serre duality}

In the following theorem we list the necessary facts concerning Serre duality. 
The first and the second statements are basic ones, see \cite{H}, \cite{C}. 
The third one follows from the very definition, see e.~g. \cite{TT}.

\begin{theorem}[Serre duality]
\begin{enumerate}
\item For $X$ proper, there is a map $\int\colon H^{\dim X}(\omega)\to k$
such that for any $E \in \D$ the composition of maps
\begin{equation}
H^i(E)\otimes H^{\dim X-i}(E^\vee\na \omega)
\stackrel{\Tr}{\longrightarrow}H^{\dim X}(\omega)\stackrel{\int}{\longrightarrow}k
\label{sd}
\end{equation}
gives a perfect (super)symmetric pairing.

\item There exists a morphism in the derived category  $\D(X\times X)$
\begin{equation}
can \colon \Oc_\Delta \to \Oc\boxtimes\omega [\dim X]
\label{can}
\end{equation}
called \emph{the canonical extension},
such that, for $X$ proper, $E,F\in \D$ and $m\in Hom(E, F)$, the composition
\[
\Oc_{X\times X}\xrightarrow{m} (F\boxtimes E^\vee)_\Delta
\xrightarrow{can \otimes (F\boxtimes E^\vee)} F\boxtimes (E^\vee\na\omega) [\dim X]
\]
is equal to
\begin{multline*}
H^*(m)\in Hom( H^*(E), H^*(F))=
H^*(F) \otimes H^*(E)^\vee\stackrel{\eqref{sd}}{=} \\H^*(F) \otimes H^{\dim X-*}(E^\vee\na\omega)=
H^*(F\boxtimes (E^\vee\na\omega)).
\end{multline*}
\item For $X= \mathop{\rm Spec} k\, [x_1, \dots, x_n]$
consider the Koszul resolution $K^\bullet$ of the diagonal \cite{BH}.
Then, the canonical extension is represented by the natural isomorphism $K_n=\Oc\boxtimes\omega$.
\end{enumerate}
\label{serre}
\end{theorem}

For $E\in \D$ let $\K$ denote the composition
\begin{equation}
 \Oc_{X\times X}\stackrel{\mathds 1}{\longrightarrow} E\boxtimes E^\vee\na\Oc_\Delta
   \stackrel{can \otimes E\boxtimes E^\vee}{\longrightarrow} E\boxtimes (E^\vee\na\omega) [\dim X]
   ,
\label{short}
\end{equation}
where the first arrow is given by the identity operator $\Oc \to
E\otimes E^\vee $. By the  second statement of the theorem,
$\K\in H^*(E)\otimes
H^*(E^\vee\otimes\omega)=H^*(E)\otimes H^*(E)^\vee=\mathop{End}
H^*(E)$ is equal to the identity operator.

By the first statement of the theorem the trace of the restriction
of $\K$ to the diagonal $\Delta^*\K\in  H^*(E\otimes E^\vee\otimes
\omega)$ followed by $\int$ is equal to the supertrace of the
identity operator on $H^*(E)$, that is, to the \emph{Euler
characteristic} :
\begin{equation*}
\chi(E)\stackrel{def}{=}\sum_i (-1)^i \dim H^i(E).
\end{equation*}

To state the Riemann-Roch theorem we need to factor the morphism $\int\Delta^*\K\in  H^*(E\otimes E^\vee\otimes \omega) $. Restricting \eqref{short} to the diagonal and taking the trace we obtain
\begin{equation}
\Oc_{X}\xrightarrow{\mathds 1\otimes\Oc_\Delta} E\na E^\vee\na\Fc
  \xrightarrow{\id\otimes (can\otimes\Oc_\Delta)} E\na E^\vee\na\omega [\dim X]
  \xrightarrow{\Tr}   \omega[\dim X],
\label{long}
\end{equation}
where $\Fc$ is defined by \eqref{tc}. Interchanging the last two arrows we obtain
\[
\Oc_{X}\xrightarrow{\mathds 1\otimes\Oc_\Delta} E\na E^\vee\na\Fc
  \xrightarrow{\Tr\otimes\id}\Fc
  \xrightarrow{can\otimes\Oc_\Delta} \omega[\dim X].
\]

We introduce the following notations: let
\begin{equation}
\Ch (E) \colon \, \Oc_{X}\xrightarrow{\mathds 1\otimes\Oc_\Delta} E\na E^\vee\na\Fc
\xrightarrow{\Tr\otimes \id}  \Fc
\label{Ch}
\end{equation}
and let
\begin{equation}
\Td \colon \, \Fc \xrightarrow{can\otimes\Oc_\Delta}  \omega [\dim X].
\label{Td}
\end{equation}

\begin{theorem}[Riemann-Roch theorem]
For $E\in \D$
\[
\chi(E)=\int \Td\circ \Ch(E)  .
\]
\end{theorem}

The classes $\Ch$ and $\Td$ are calculated explicitly in propositions
\ref{ch} and \ref{td} below.

\subsection{The Chern character}

For $E$ in $\D$ let
\begin{equation}
\at\nolimits^i(E)\colon\quad E\longrightarrow E\na {(\Omega^1)}^{\otimes i}[i]
\xrightarrow{\id\otimes\wedge} E \na \Omega^i [i],
\label{ae}
\end{equation}
where the first arrow is the $i$-fold of $\at(E)$ with itself, and the
second one is the usual multiplication in $\Omega^*$. Let
\[
\wedge \at\nolimits^i(E)  \colon\quad   E \na \Omega^*
\xrightarrow{\at\nolimits^i(E)\otimes\id}
 E \na \Omega^i \otimes  \Omega^*[i]
\xrightarrow{\id\otimes\wedge} E \na \Omega^{i+*}[i].
\]  
We define
\[
\exp(\at(E)) \in \bigoplus_i Ext^i(E, E\na \Omega^i[i])
\]
by the formula  $\exp(\at(E))=\sum \at^i(E)/i!$ and
\[
\wedge \exp(\at(E)) \colon\quad E\na \bigoplus_i \Omega^i[i]
\longrightarrow    E\na \bigoplus_i \Omega^i[i].
\]
by the formula  $\wedge\exp(\at(E))=\wedge\sum \at^i(E)/i!$.

For $i = 1,2$ we define the isomorphism $a_i$ as the composition of isomorphisms
\begin{multline*}
 a_i\colon E_\Delta\na \Oc_\Delta = (p_i^*E\na
\Oc_\Delta)\na\Oc_\Delta= \\
p_i^*E\na(\Oc_\Delta\na\Oc_\Delta)=p_i^*E\na\Fc_\Delta=(E\na \Fc)_\Delta.
\end{multline*}
Let
\[
\E E_i\colon Rp_{1*}(E_\Delta\na
\Oc_\Delta)\stackrel{a_i}{=}E\na\Fc\stackrel{\id\otimes\E}{=}E\na\bigoplus_j
\Omega^j[j]
\]
denote the composition of the image of $a_i$ under $Rp_{1*}$ with $\E$.

The following lemma could be an alternative definition of the Atiyah
class.

\begin{lemma}
In the notation introduced above we have:
\begin{equation}
\E E_1= (\id\otimes \wedge \exp(\at(E)))\circ \E E_2. \label{eq}
\end{equation}
\label{12}
\end{lemma}

\begin{proof} By analogy with  \eqref{ac} we define an action of
$\Tc$ on $Rp_{1*}(E_\Delta\na\Oc_\Delta)$
\begin{equation}
\Tc\na Rp_{1*}(E_\Delta\na\Oc_\Delta)\rightarrow Rp_{1*}(E_\Delta\na\Oc_\Delta)
\label{198}
\end{equation}
as the action of  $(\id\otimes\pi)\circ \at(E_\Delta)\in
Ext^1(E_\Delta, E_\Delta\na p_1^* \Omega^1)$ on $E_\Delta$
(where $\pi\colon
\Omega^1_{X\times X} \to p_1^* \Omega^1_X $ is the projection) followed by the restriction to the diagonal.

Substituting \eqref{198} into the definition of the isomorphisms $a_i$ and
using the definition of the action of $\Tc$ on $\Fc$ in \eqref{actt} and
Lemma \ref{na} one obtains the commutative diagrams
\begin{equation}
\begin{CD}
\Tc\na Rp_{1*}(E_\Delta\na\Oc_\Delta) @>\eqref{198}>> Rp_{1*}(E_\Delta\na\Oc_\Delta)\\
@V\id\otimes Rp_{1*} (a_1)VV @VRp_{1*}(a_1)VV \\
 \Tc\na (E\na\Fc)   @>\Dc>>E\na\Fc
\end{CD}
\label{com1}
\end{equation}
 where $\Dc$ denotes the morphism \eqref{actt}, and
\begin{equation}
\begin{CD}
\Tc\na Rp_{1*}(E_\Delta\na\Oc_\Delta) @>\eqref{198}>> Rp_{1*}(E_\Delta\na\Oc_\Delta) \\
@V\id\otimes Rp_{1*}(a_2)VV @VRp_{1*}(a_2)VV \\
 \Tc\na (E\na\Fc)  @>\Dc+\at(E)>>E\na\Fc
\end{CD}
\label{com2}
\end{equation}
The latter diagram is commutative because the projection of $\at(p^*_i E)$ to $p^*_i\Omega^1$ is equal to $p^*_i \at(E)$ while the projection to the complementary bundle is zero.

Arguments as in Theorem \ref{un} show that action of $\Tc$ \eqref{198}
extends to an action of $\Uc$:
\begin{equation}
\Uc \na Rp^*_1( E_\Delta\na\Oc_\Delta) \rightarrow Rp^*_1( E_\Delta\na\Oc_\Delta) .
\label{*}
\end{equation}
Let 
\[
\varepsilon_E \colon Rp^*_1( E_\Delta\na\Oc_\Delta) \rightarrow E
\]
denote the morphism induced by the canonical morphism $E_\Delta\na\Oc_\Delta \to E_\Delta$.
Combining  \eqref{*} with $\varepsilon_E $ one obtains the morphism
\[
\Uc \na Rp^*_1( E_\Delta\na\Oc_\Delta) \rightarrow E
\]
and, dually, the morphism
\begin{equation}
Rp^*_1( E_\Delta\na\Oc_\Delta) \rightarrow E\na\Uc^\vee\stackrel{\id\otimes\E}{=}E\na\bigoplus_i
\Omega^i[i].
\label{fe}
\end{equation}

It follows from the definition of $\E$ and \eqref{com1} that the composition \eqref{fe} is equal to the left hand side of \eqref{eq}. It follows from \eqref{com2} that it is equal to the right
side of \eqref{eq}. This proves the lemma. \end{proof}

\begin{definition} For  $E\in \D$, the \emph{Chern character} is defined by the formula
\[
\ch(E)=\Tr \sum_i \at\nolimits^i(E)/i! \in \bigoplus_i H^i(\Omega^i).
\]
\label{chern}
\end{definition}

\begin{prop}
For $E\in \D$ the composition of the map  $\Ch(E)$ given by \eqref{Ch} and the isomorphism $\E$
is equal to $\ch(E)$.
\label{ch}
\end{prop}
\begin{proof}
It is clear from Definition \ref{chern} that it is sufficient to show that the composition of the first arrow in \eqref{Ch} with the  map $\E$  is equal to $\exp(\at(E))$.

For $i,j = 1,2$ we have isomorphisms $I_{i,j}$ defined as the compositions
\[
I_{i,j}\!:\, Rp_{1*}((E\na E^\vee)_\Delta\na \Oc_\Delta) =
Rp_{1*}( (p_i^*E\na p_j^*E^\vee)\na\Oc_\Delta) =
(E\na E^\vee)\na \Fc,
\]
The map in question is the composition of the identity section
$\mathds1\colon \Oc_{X\times X} \to (E\na E^\vee)_\Delta   $
restricted to the diagonal with the isomorphism $I_{1,2}$:
\begin{equation}
 \Oc_X \xrightarrow{\mathds 1\otimes\Oc_\Delta} Rp_{1*}((E\na E^\vee)_\Delta\na \Oc_\Delta)
 \xrightarrow{I_{1,2}} E\na E^\vee \na \Fc  .
\label{l1}
\end{equation}

The identity section of $(E\na E^\vee)_\Delta$ is equal to the composition
\begin{equation}
\Oc_{X\times X}\xrightarrow{1} \Oc_\Delta
\xrightarrow{p^*_1 \mathds 1\otimes \Oc_\Delta} p_1^*(E\na E^\vee)\na \Oc_\Delta
=(E\na E^\vee)_\Delta,
\label{i1}
\end{equation}
where $\mathds 1\colon \Oc_X \to E\na E^\vee $ is the identity and
$1\colon \Oc_{X\times X}\to \Oc_\Delta$ is the canonical map.
Restricting to the diagonal and applying $Rp_{1*}$ one obtains
\begin{equation}
\Oc_X \xrightarrow{1\otimes\Oc_\Delta} \Fc
\xrightarrow{\mathds 1\otimes \id} E\na E^\vee\na \Fc,
\label{i2}
\end{equation}
using $I_{1,1} \colon Rp_{1*}((E\na E^\vee)_\Delta\na \Oc_\Delta) = (E\na E^\vee)\na \Fc$.
Thus, the composition
\begin{equation}
\Oc_X \stackrel{\mathds 1\otimes\Oc_\Delta}{\longrightarrow} Rp_{1*}((E\na E^\vee)_\Delta\na \Oc_\Delta)
 \stackrel{I_{1,1}}{\longrightarrow } E\na E^\vee \na \Fc  .
\label{l2}
\end{equation}
is equal to $\mathds 1\otimes 1$, where $1\colon \Oc\to \Fc$ is the natural embedding.

Applying Lemma \ref{12} to \eqref{l1} and \eqref{l2} one proves the statement.

\end{proof}

\subsection{The Todd class}

\begin{prop} Let
\begin{equation}
\td=\exp(\sum t_i \ch(\Omega^1))\in \bigoplus_iH^i(\Omega^i),
\label{t}
\end{equation}
where $\sum t_iz^i= \log(z/(e^z-1))$.
Then, the class $\Td$ from \eqref{Td} may be expressed as the composition
\[
\Td\colon \Fc \stackrel{\E}{\longrightarrow} \bigoplus_i\Omega^i[i] \stackrel{\wedge \td}{\longrightarrow}
\bigoplus_i\Omega^i[i] \twoheadrightarrow \omega[\dim X],
\]
where the last arrow is the projection onto the differential forms of top degree.
\label{td}
\end{prop}
\begin{proof} Applying \eqref{com} to the canonical map \eqref{can} and composing with the projection
$\pi\colon\Omega^1_{X\times X}\to p^*_1 \Omega^1_X $ one obtains the commutative diagram
\begin{equation}
\begin{CD}
\Oc_\Delta                        @>can>>            \Oc\boxtimes\omega [\dim X]\\
@V\pi\circ\at(\Oc_\Delta) VV                                   @V\at(\Oc)VV \\
\Oc_\Delta\otimes p_1^* \Omega^1 @>can\otimes\id>> \Omega^1\boxtimes\omega [\dim X+1]
\end{CD}
\end{equation}
The right vertical arrow is zero, using Lemma \ref{na} and the fact that the Atiyah class of $\Oc$ is
trivial. Therefore, so is the composition of the bottom and the  left vertical arrows. Tensoring by $\Oc_\Delta$ one finds (using notation
of \eqref{actt}) that the composition
\begin{equation}
\Td\circ\Dc\colon \,\Fc \xrightarrow{\Dc} \Fc\otimes \Omega^1 [1]
\xrightarrow{\Td\otimes \id}  \omega\otimes \Omega^1 [\dim X+1]
\label{zero}
\end{equation}
is equal to zero.

Applying  isomorphism $\E$ to \eqref{zero}  and making use of notations  of Theorem \ref{L}
we obtain
\[
\bigoplus_i\Omega^i[i] \xrightarrow{\Lc}
(\bigoplus_i\Omega^i[i]) \otimes \Omega^1 [1] \xrightarrow{x\otimes\id}\omega\otimes\Omega^1[\dim X+1].
\]
It follows from Lemma \ref{l} below that vanishing of \eqref{zero} determines $\Td$
up to a scalar factor. This factor (which is 1) can be determined
from the local considerations by means of the third part of Theorem \ref{serre}.
This finishes the proof.\end{proof}

\begin{lemma} Suppose that  $x\colon\bigoplus_i\Omega^i[i]\to \omega[\dim X]$
is a morphism such that the composition
\[
\bigoplus_i\Omega^i[i] \xrightarrow{\Lc}
(\bigoplus_i\Omega^i[i]) \otimes \Omega^1 [1] \xrightarrow{x\otimes\id}\omega\otimes\Omega^1[\dim X+1]
\]
is zero. Then, up to a factor it is given by the composition of $\wedge\td$ (see \eqref{t}) with
the projection onto the differential forms of top degree.
\label{l}
\end{lemma}

\begin{proof}
Denote by $\overline{\,\cdot\,}$ the anti-involution on
$\bigoplus_i\Omega^i[i]$ which multiplies $\Omega^i$ by $(-1)^{i(i-1)/2}$.
We define a non-degenerate pairing on $\bigoplus_i\Omega^i[i]$ by the formula
\begin{equation}
\langle,\rangle\colon(\bigoplus_i\Omega^i[i])\otimes(\bigoplus_i\Omega^i[i])
\stackrel{\overline{\,\mathstrut\cdot\,}\wedge\,\mathstrut\cdot\,}{\longrightarrow}
\bigoplus_i\Omega^i[i] \twoheadrightarrow \omega[\dim X]
\end{equation}
where the last arrow is the projection onto the forms of top degree.
There exists a unique $y\in \bigoplus_iH^i(\Omega^i)$ such that $\langle \overline{y}, \,\mathstrut\cdot\,\rangle=x$.

Let $\Lc^+$ denote the map adjoint to $\Lc$ with respect to
$\langle,\rangle$. This means that the diagram
\begin{equation}
\begin{CD}
(\bigoplus_i\Omega^i[i])\otimes(\bigoplus_i\Omega^i[i]) @>\id\otimes \Lc>> (\bigoplus_i\Omega^i[i])\otimes(\bigoplus_i\Omega^i[i])\otimes\Omega^1[1]\\
@V\Lc^+\otimes \id VV    @V\langle,\rangle \otimes\id VV \\
(\bigoplus_i\Omega^i[i])\otimes(\bigoplus_i\Omega^i[i])\otimes\Omega^1[1] @>\langle,\rangle\otimes\id>>
\omega\,[\dim X]\otimes\Omega^1[1]
\end{CD}
\end{equation}
is commutative.
It follows from Theorem \ref{L} by direct calculation that
$\Lc^+=\wedge\overline{\td}\circ\Lc\circ\wedge\overline{\td}^{-1} $.
(That is, $\td$ is analogous to the left invariant volume form on a Lie group.)

By the hypothesis of the lemma $0=\langle \overline{y}, \Lc\,\mathstrut\cdot\,\rangle=\langle\Lc^+
\overline{y},\,\mathstrut\cdot\,\rangle=
\langle\overline{\td}\wedge\Lc\overline{\td\nolimits^{-1}\wedge y}, \,\mathstrut\cdot\,\rangle$.
Since the pairing is non-degenerate and $\overline{\td}$ is invertible, it follows that
$\Lc(\td^{-1}\wedge y)=0$.

Thus, it remains to prove that the only sections $s\colon\Oc\to\bigoplus_i\Omega^i[i]$
for which the composition
\begin{equation}
\Oc\xrightarrow{s}\bigoplus_i\Omega^i[i]\xrightarrow{\Lc}
\bigoplus_i\Omega^i[i]\otimes\Omega^1[1]
\label{sec}
\end{equation}
vanishes are  the ones which
factor through $\Oc\hookrightarrow\bigoplus_i\Omega^i[i]$.
It follows from formula \eqref{formula} that the component of $\Lc$
in $Hom(\Omega^i[i], \Omega^{i-1}[i-1]\otimes\Omega^1[1]) $ for $i>0$ is equal to the embedding
$\Omega^i\hookrightarrow\Omega^{i-1}\otimes \Omega^1 $ as the skew-symmetric part.
This means that \eqref{sec} cannot vanish if $s$ has components
in $\Omega^{>0}$.
\end{proof}

\bibliographystyle{abbrv}
\bibliography{markarian}

\end{document}